\documentclass[12pt]{article}
\usepackage{amsmath,amssymb,amsthm,amscd,latexsym}
\usepackage{epsfig,subfigure,graphicx}
\newcommand{\pea}{{\mathbf p}}
\newcommand{\ex}{{\mathbf x}}

\newtheorem{thm}{Theorem}[section]

\newtheorem{Def}[thm]{Definition}

\begin{document}

\title{Thermoacoustic Tomography - Implementation of Exact Backprojection Formulas}
\author{Gaik Ambartsoumian\dag\ and Sarah K Patch\ddag\\}
\date{}
\maketitle

\begin{abstract}
The problem of image reconstruction in thermoacoustic tomography requires inversion of a generalized Radon transform, which integrates the unknown function over circles in 2D or spheres in 3D. The paper investigates implementation of the recently discovered backprojection type inversion formulas for the case of spherical acquisition in 3D. A numerical simulation of the data acquisition with subsequent reconstructions are made for the Defrise phantom as well as for some other phantoms. Both full and partial scan situations are considered. The results are compared with the implementation of a previously used approximate inversion formula.
\end{abstract}


\section{Introduction}

\begin{figure}
\begin{center}
{\epsfig{figure= 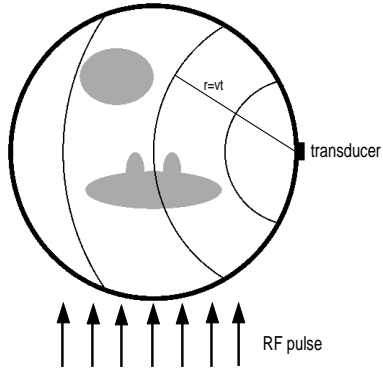,width=2in, height=2in }}
\end{center}
\caption{Sketch of a TCT system}\label{fig:TCTsketch}
\end{figure}

The idea of thermoacoustic tomography (TCT, sometimes also called TAT)~\cite{ref:KrugerRadiologyMedPhys99, ref:KrugerRadiology99, ref:KrugerRadiology00, ref:Wang, ref:XWAK} can be briefly described as follows (see Figure~\ref{fig:TCTsketch}). A short pulse of radiofrequency (RF) electromagnetic waves is sent through a biological object heating up the tissue. It is known that the cancerous cells absorb several times more RF energy than the healthy ones~\cite{ref:Joines94}. As a result a significant increase of temperature occurs at the tumor locations causing a thermal expansion of cancerous masses pressing on the neighboring healthy tissue. The created pressure wave is registered by the transducers located on the edge of the object. Assuming speed of propagation of these acoustic waves constant inside the object (an assumption which is not always correct, but is satisfactory, e.g. for mammography), the signals registered at any transducer location are generated by the inclusions laying on a sphere centered at that location. In fact the measured data are the integrals of RF absorption coefficient $f$ over those spheres or, in other words, the spherical Radon transform $Rf$ of the RF absorption coefficient $f$. Hence, to reconstruct the image one needs to invert the spherical Radon transform.

\begin{Def}
The spherical Radon transform of $f$ is defined as
$$
Rf(p,r)=\int_{|x-p|=r}f(x)d\sigma(x),
$$
where $d \sigma(x)$ is the surface area on  the sphere $|x-p|=r$ centered at $p\in \mathbb{R}^3$.
\end{Def}

In the definition we allow arbitrary set of centers $p$ and radii
$r$. However, from the dimensional consideration it is clear, that
this mapping is overdetermined. The tomographic motivation suggests
to restrict the set of centers to a surface $S \subset \Bbb{R}^3$
(the set of transducers' locations), while not imposing any
restrictions on the radii. In this paper we will deal only with the
case of spherical acquisition (i.e. the transducers are located on a
unit sphere) and from now on we will suppose $|p|=1$.

Two different approaches have been used to derive exact inversion
formulae for this case.  Fourier-Bessel and spherical harmonic
expansions result in solutions written as an infinite series for two
and three dimensions respectively~\cite{ref:Norton2D,ref:Norton3D}.
The TCT analog of $\rho$-filtered backprojection inversion is
derived in~\cite{ref:TCTSIAM}
\begin{equation}
f(x)=-\frac{1}{8\pi^2}\triangle_x\left( \int_{|p|=1} \frac{1}{|x-p|}\; Rf(p,|x-p|)dp \right)  \label{eq:rhofilt}
\end{equation}
as well as a filtered backprojection type version
\begin{equation}
f(x)=-\frac{1}{8\pi^2}\left( \int_{|p|=1} \frac{1}{|x-p|}\;\; \frac{\partial^2}{\partial r^2} Rf(p,|x-p|)dp \right)  \label{eq:FBP}
\end{equation}
Both formulas can be generalized to higher odd
dimensions~\cite{ref:TCTSIAM}. Notice that, as one can expect for a
codimension 1 Radon transform in 3D, the formulas are local.

In Section 2 we describe the numerical simulation of the data
acquisition. The reconstruction algorithms based on the
$\rho$-filtered backprojection formula~(\ref{eq:rhofilt}) and the
filtered backprojection one~(\ref{eq:FBP}) are discussed in
Section~3.


\section{Data Simulation} \label{sec:forward}

\begin{figure}[h]
\begin{center}
{\epsfig{figure=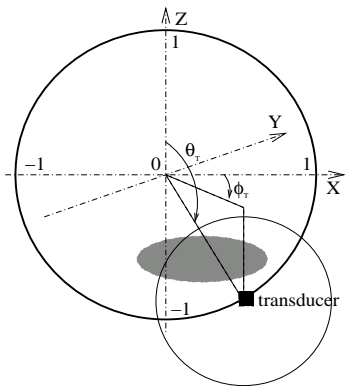,width=1.8in,height=2in }}
\end{center}
\vspace{-0.7cm}
\caption{Phantom setup}\label{fig1}
\end{figure}

The region of reconstruction is the unit ball centered at the origin (see fig.~\ref{fig1}). All phantoms considered in the paper are sums of indicator functions supported in ellipsoids completely contained inside the unit ball. The transducers are located on the surface of the unit sphere.

We parameterize the transducer location by two angles $(\phi_T,\theta_T)$, where $\phi_T\in [0,2\pi)$ is the the azimuthal angle in the $xy$-plane and $\theta_T\in[0,\pi]$ is the polar angle measured from the $z$-axis.

Since the spherical Radon transform is linear, it is enough to create projections for phantoms with a single ellipsoid and then superimpose the projections. For a single ellipsoid the data measured at a fixed transducer location at a given moment (i.e. for fixed $(\phi_T,\theta_T,r)$) is the surface area of a part of the sphere of integration cut by the intersecting ellipsoid. It can be expressed as a finite sum with terms of the form
\begin{equation}
\int_0^{2\pi}\int_{{\theta_1}({\phi})}^{{\theta_2}({\phi})} \sin{\theta} \; d{\theta}\; d{\phi} = \int_0^{2\pi} \left[ \cos{\theta_1}({\phi}) - \cos{\theta_2}({\phi}) \right] d{\phi} \label{eq:doubleint}
\end{equation}
where each such term corresponds to a connected component of the intersection. Here ${\phi}$ and ${\theta}$ parameterize the sphere of integration and are independent of $\phi_T$ and $\theta_T$, which parameterize the transducer location. The angles ${\theta_1}({\phi})$ and ${\theta_2}({\phi})$ are defined by the intersection of the integration sphere and the phantom's ellipsoid. The cosines of these angles can be found from the solution of a quartic equation describing that intersection.

\begin{figure}[h]
\begin{center}
{\epsfig{figure= 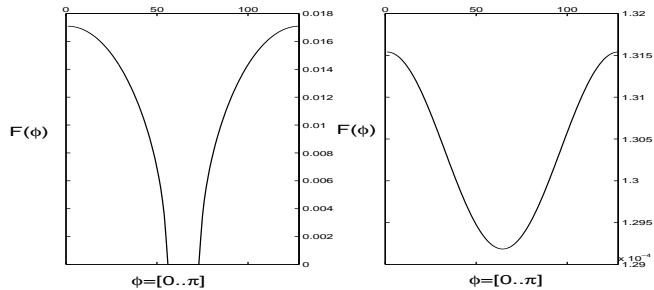,width=3.4in,height=1.5in }}
\end{center}
\vspace{-.5cm}
\caption{$F(\phi)$ for two different locations of an ellipsoid}
\label{fig4}
\end{figure}

In the numerical results presented below the quartic equation is solved using the MATLAB built-in function ``roots".  By adding up these roots in an appropriate way we obtain the inner integral with respect to the polar angle $\theta$ in equation~(\ref{eq:doubleint}). The result is a function of azimuthal angle $\phi$, which we will denote $F(\phi)$. Depending on the location and parameters of the ellipsoid, $F(\phi)$ might be either a smooth $\pi$-periodic function of $\phi$, or a piecewise smooth one (see fig.~\ref{fig4}).
 In the first case we compute its values at uniformly discretized locations on the interval $[0,\pi]$ and use the trapezoidal rule to compute the integral. For \mbox{$F(\phi) \in C^2$}, numerical integration using the trapezoidal rule is accurate to $O(h^4)$~\cite{ref:Atkinson}. If, however, $F(\phi)$ is only piecewise smooth on $[0,\pi]$, then we locate the pieces of \mbox{$supp\; F(\phi)$} where it is smooth and use Gaussian quadrature to integrate over each piece.


\section{Reconstruction}  \label{sec:inverse}

Once we have generated the projection data, we reconstruct the
original indicator functions of the phantoms. The reconstruction
algorithms are based on the $\rho$-filtered
backprojection~(\ref{eq:rhofilt}) or the filtered
backprojection~(\ref{eq:FBP}).

The integrals over the unit sphere in~(\ref{eq:rhofilt})
and~(\ref{eq:FBP}) are computed as double integrals with respect to
the azimuthal angle $\phi_T$ and the polar angle $\theta_T$.  The
function to be integrated is periodic with respect to $\phi_T$
making the trapezoidal rule an appealing quadrature choice.
Integration with respect to $\theta_T$ is done by Gaussian
quadrature. The Laplace operator is implemented through the Matlab
built-in function ``del2".  The reconstructions were generated using
Matlab 5.0.

In the results below the resolution is
\mbox{$256\times256\times256$} over a \mbox{$2 \times 2 \times 2$}
volume resulting in isotropic pixel dimension of \mbox{$1/128$}.

\begin{figure}[h]
\begin{center}
{\epsfig{figure= 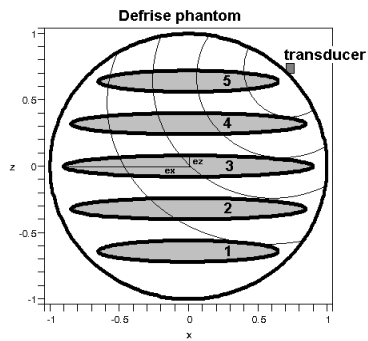,width=2.7in,height=3.5in }}
\end{center}
\vspace{-4.2cm}
\caption{The Defrise phantom slice along the plane y=0}\label{fig:Def2d}
\end{figure}

The algorithm is tested on the Defrise phantom which consists of five thin ellipsoids symmetrically centered along the $z$-axis (see fig.~\ref{fig:Def2d}). We numerate them from $1$ to $5$ starting with the lowest.

\vspace{0.5cm}

\begin{tabular}{|c|c|c|}
\hline
ellipse number & center $=(x_0,y_0,z_0)$ & semiaxes lengths $=(e_x,e_y,e_z)$\\
\hline
1 &  $(0,0,-0.64)$ & $(0.65,0.65,0.08)$\\
\hline
2 &  $(0,0,-0.32)$ & $(0.85,0.85,0.08)$\\
\hline
3 &  $(0,0,0)$ & $(0.9,0.9,0.08)$\\
\hline
4 &  $(0,0,0.32)$ & $(0.85,0.85,0.08)$\\
\hline
5 &  $(0,0,0.64)$ & $(0.65,0.65,0.08)$\\
\hline
\end{tabular}

\subsection{Full scan data}

The data was acquired from the transducers located discretely over the sphere in the following way. The azimuthal angles of the transducer locations were uniformly discretized to $N_{\phi}=400$ points between 0 and $2\pi$. The polar angles of the transducer locations corresponded to $N_{\theta}=200$ Gaussian nodes on the interval form 0 to $\pi$, as described in the previous section. The radii of the integration spheres were uniformly discretized to $N_{r}=200$ points from 0 to 2. The reconstruction was done by both methods: filtered backprojection (FBP) and $\rho$-filtered backprojection.
\begin{figure}[h]
\begin{center}
\subfigure[FBP]{\epsfig{figure= 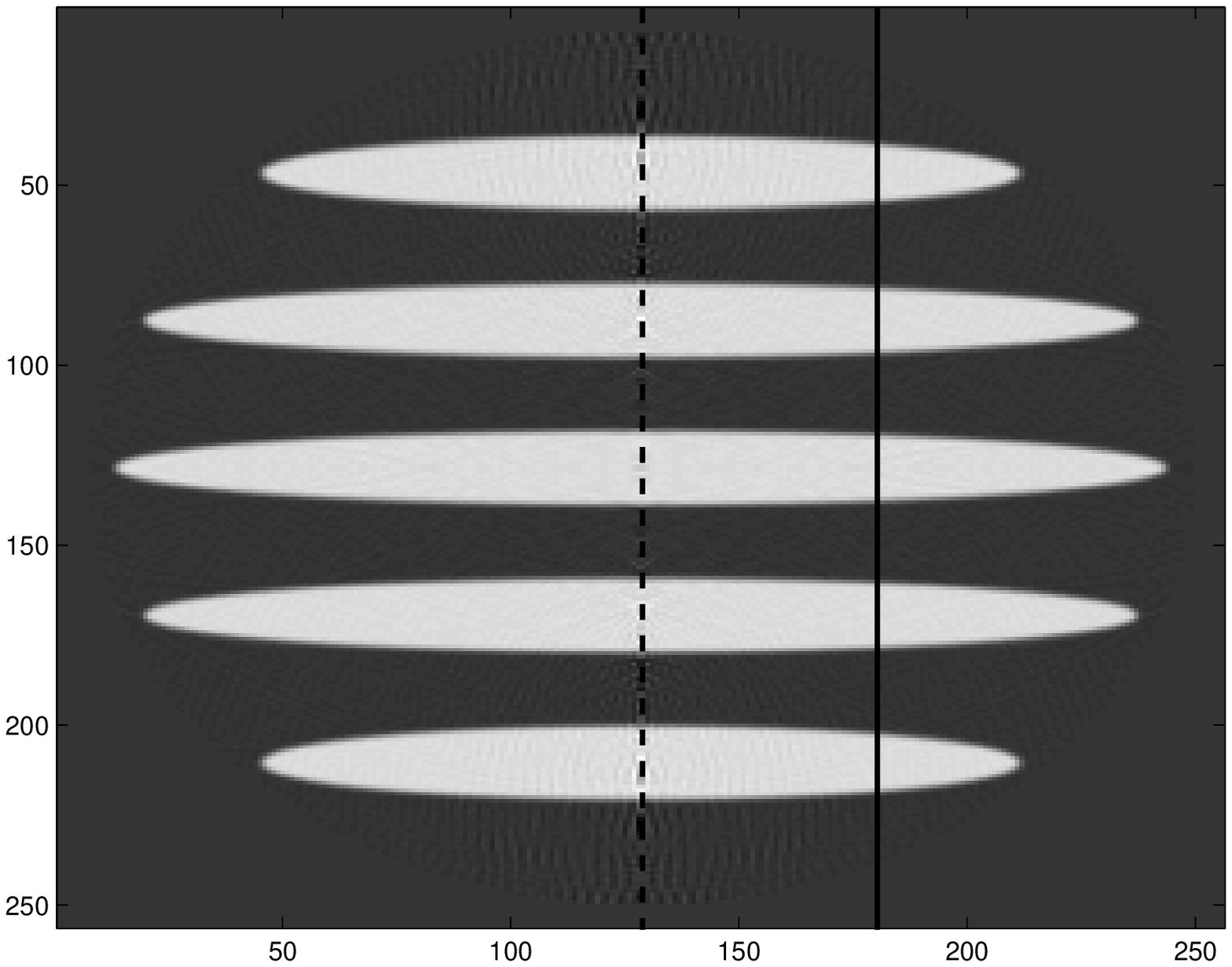,width=1.8in,height=1.8in }}
\subfigure[$\rho$-filtered reconstruction. ]{\epsfig{figure= 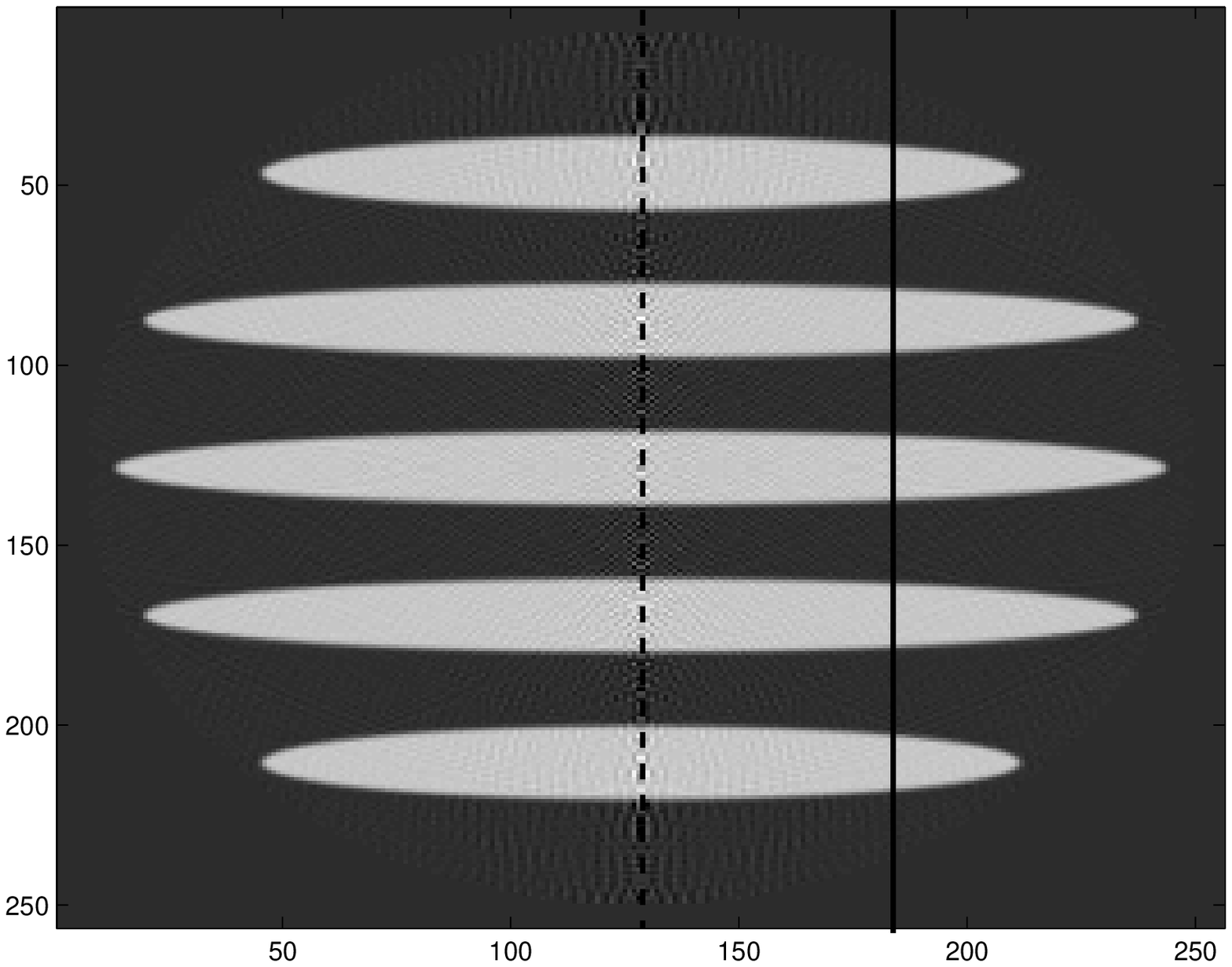,width=1.8in,height=1.8in }} \\
{\epsfig{figure= 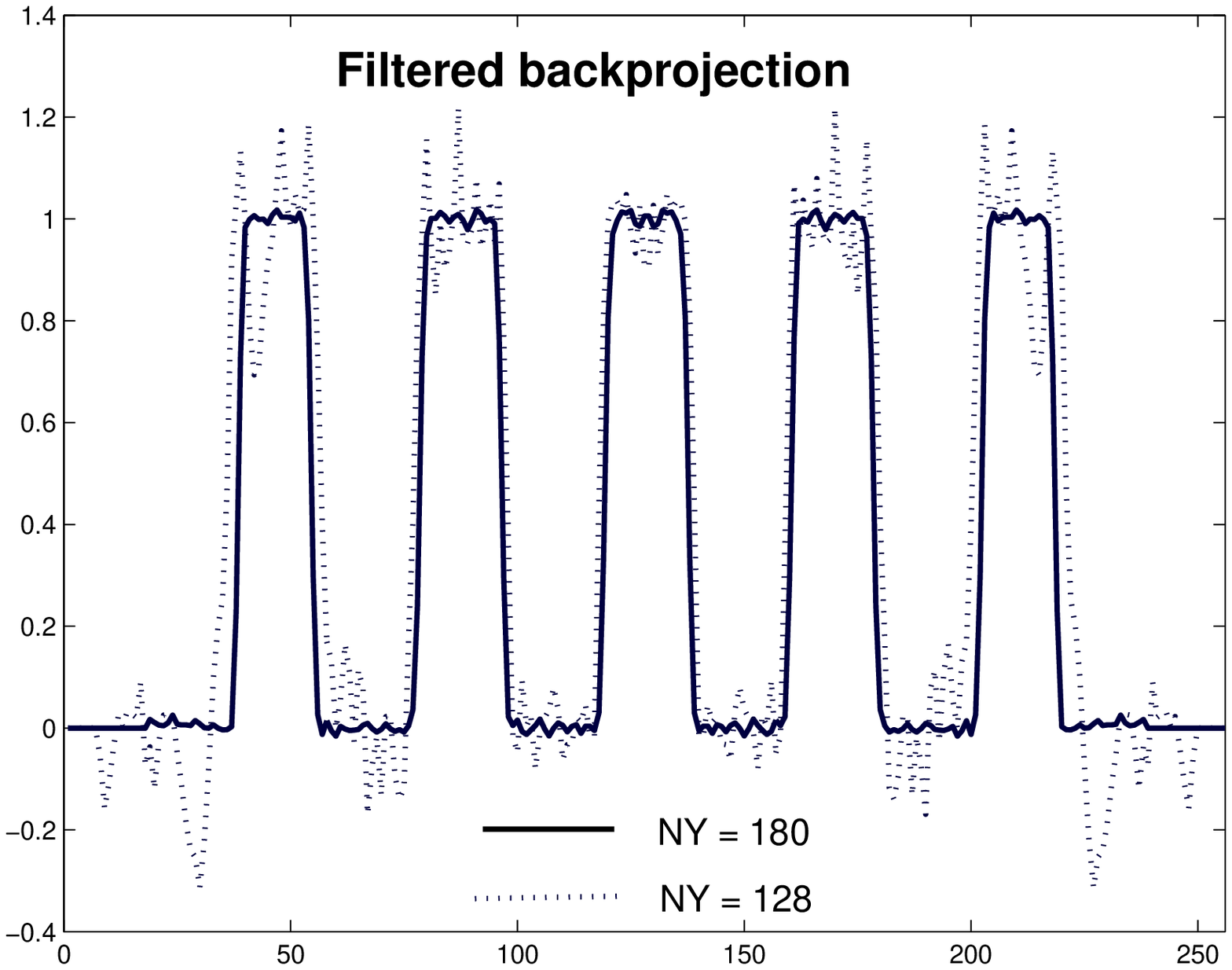,width=1.8in,height=1.8in }}
{\epsfig{figure= 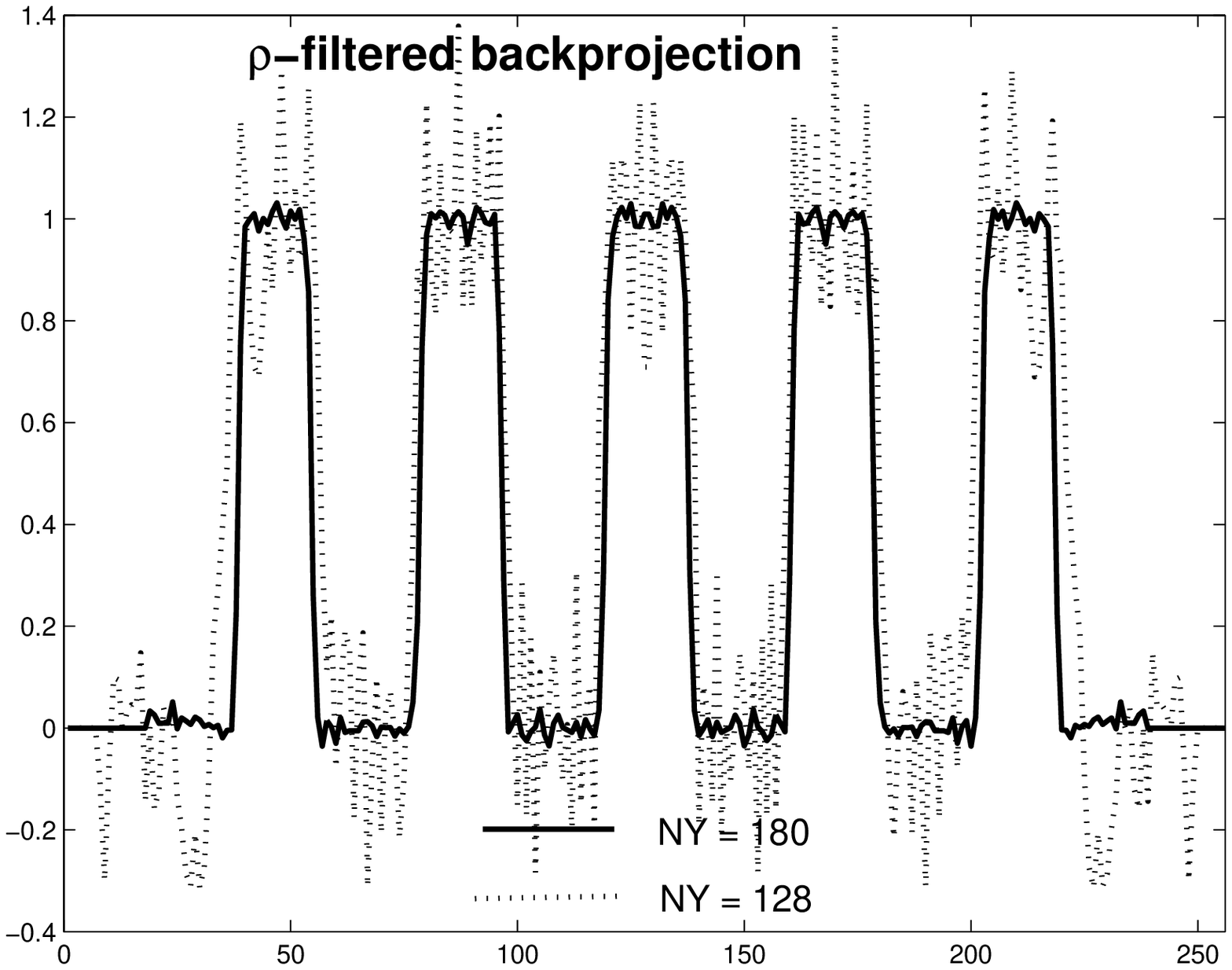,width=1.8in,height=1.8in }}
\end{center}
\caption{Reconstructions and profiles of the Defrise phantom along the center \mbox{$x=0$} slice.  Dashed lines correspond to the center \mbox{$x = 0 = y$} profile; solid lines correspond to \mbox{$x = 0, y=0.4$}} \label{Defrise}
\end{figure}

The obtained results validate reconstruction formulas~(\ref{eq:rhofilt}) and~(\ref{eq:FBP}) (see fig.~\ref{Defrise}). In both cases the Defrise phantom has a good reconstruction everywhere except along the $z$-axis ($x=y=0$), where some noise is present which, while not always noticeable on reconstructions, is visible on the graphs. The reason for appearance of that noise is the correlation of numerical errors along that axis of phantom's symmetry and is discussed in section~\ref{subsec:errors}.

\subsection{Partial scan data}
\begin{figure}[h]
\begin{center}
\subfigure[$\rho$-filtered BP; $1/2$-scan with respect to $\phi$.  The profiles along the lines are plotted below. ]{\epsfig{figure= 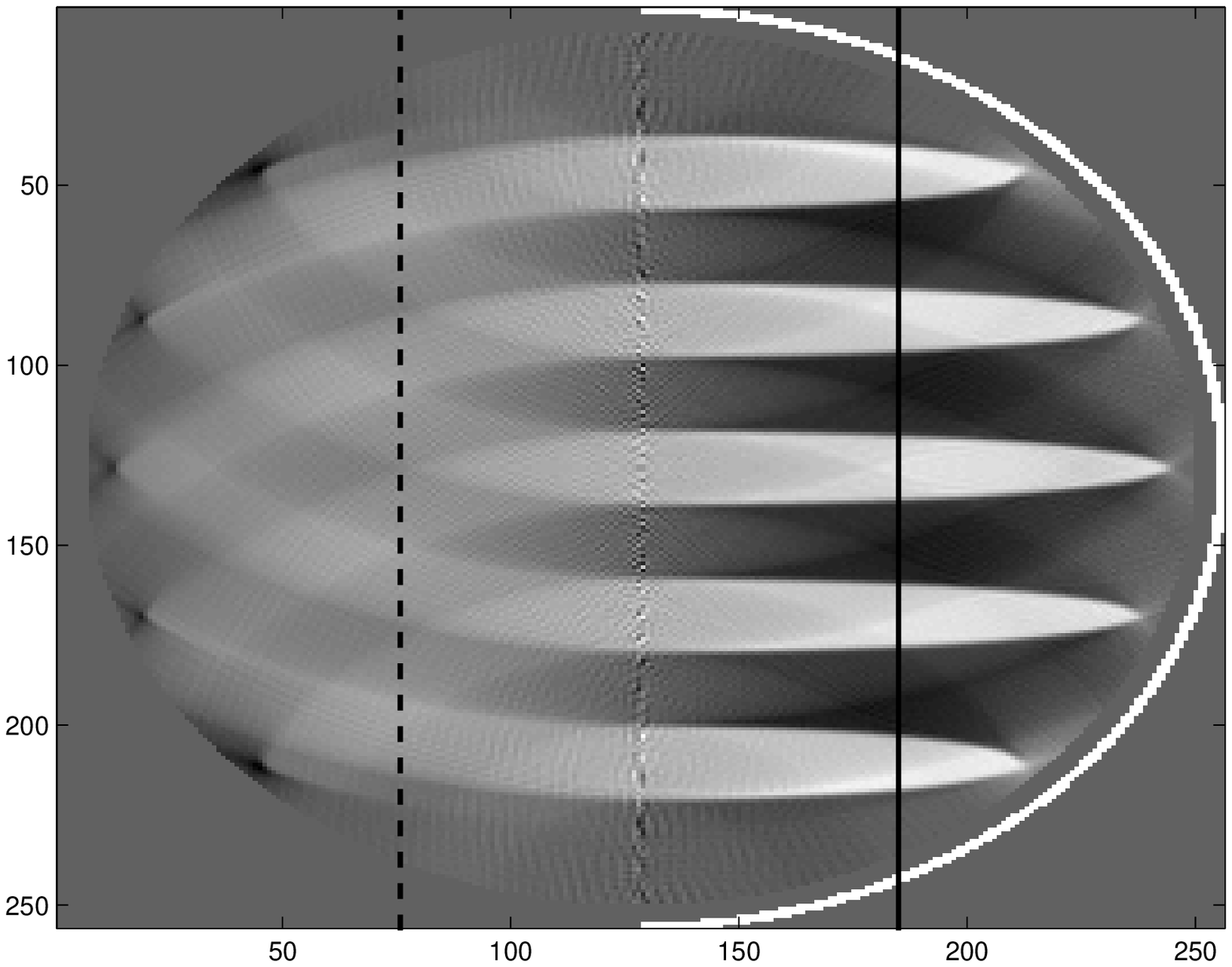,width=1.7in,height=1.7in }}
\subfigure[FBP; $1/2$-scan with respect to $\theta$.  The profile along the solid line is plotted below.]{\epsfig{figure= 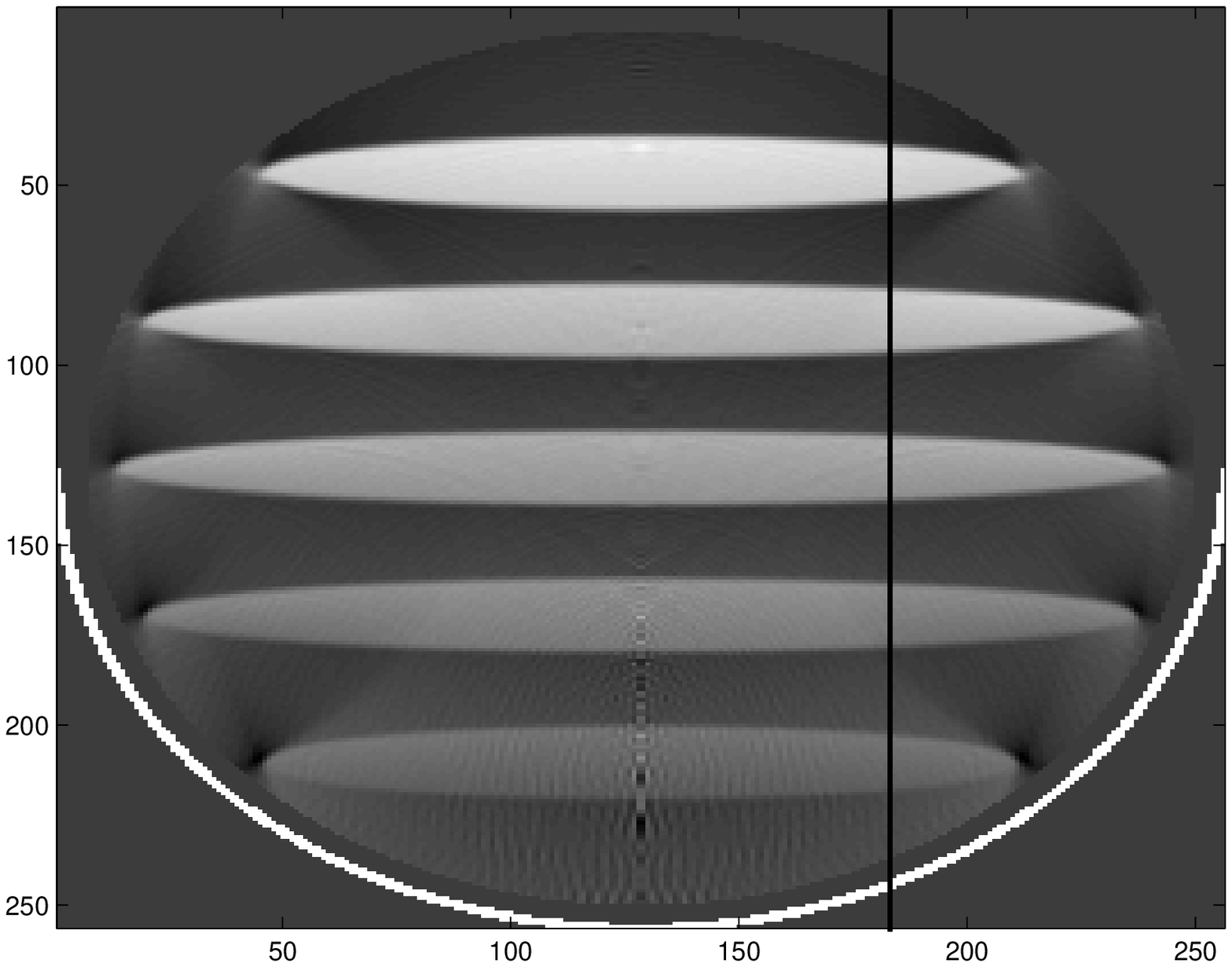,width=1.7in,height=1.7in }} \\
{\epsfig{figure= 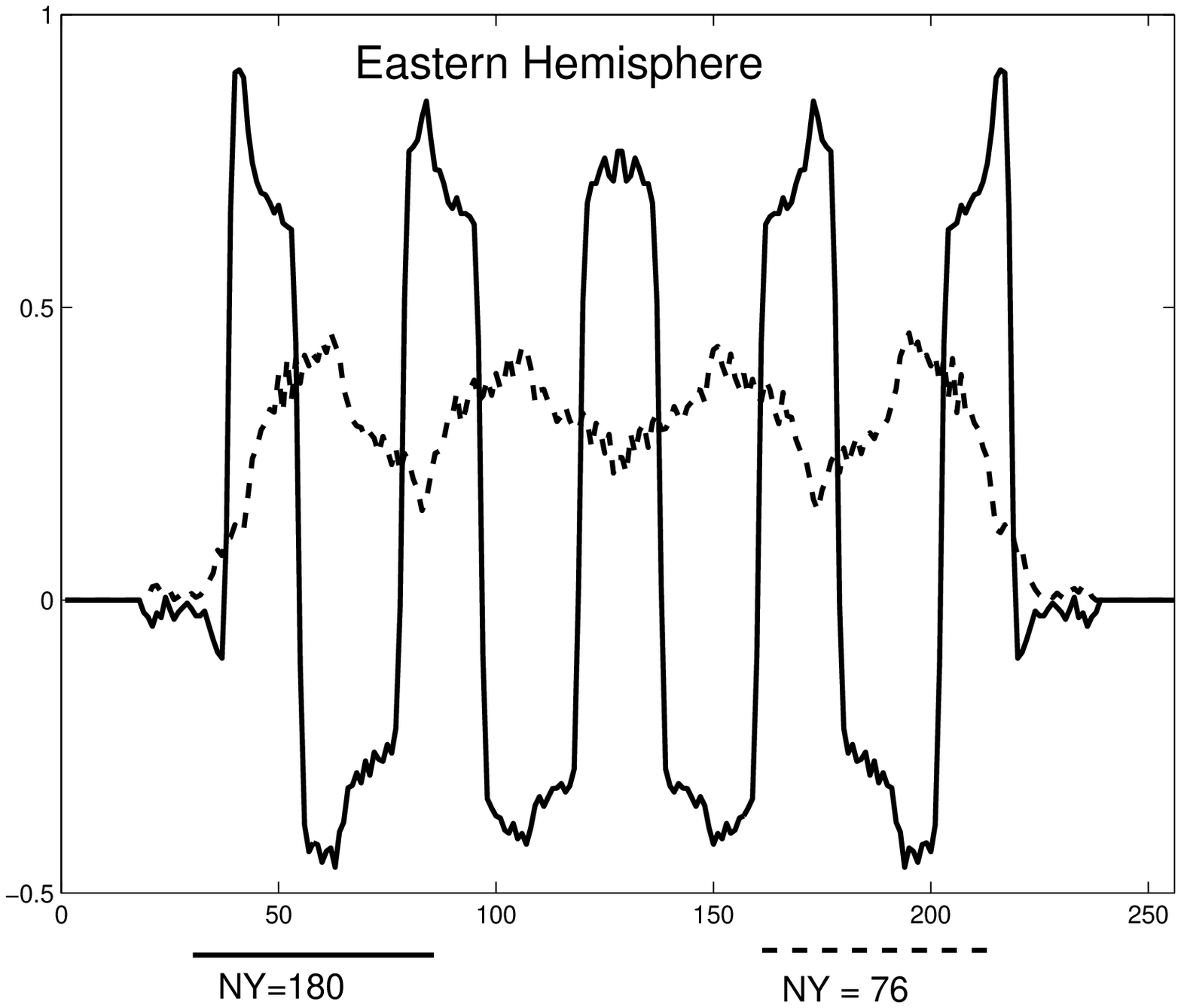,width=1.7in,height=1.7in }}
{\epsfig{figure= 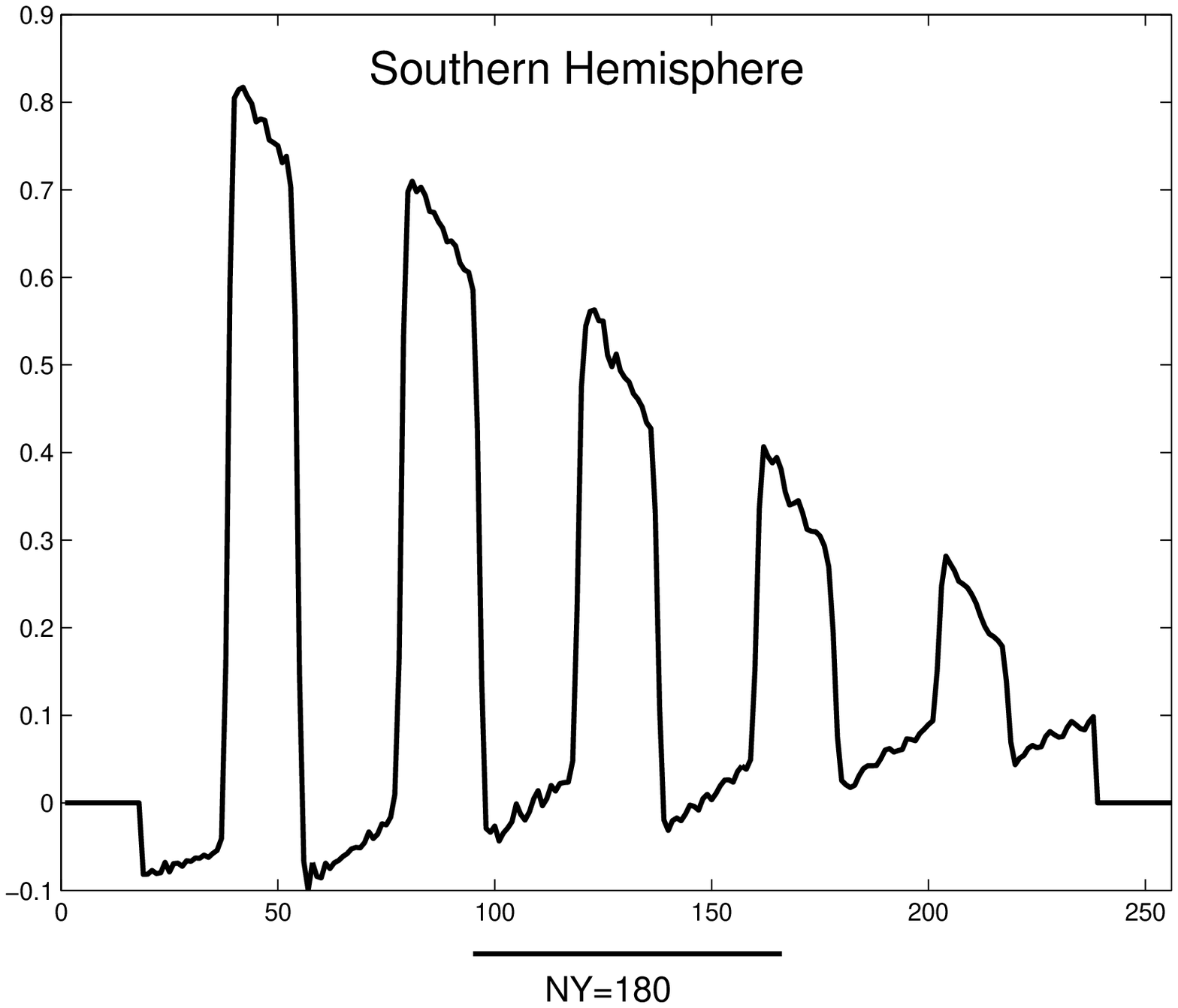,width=1.7in,height=1.7in }}
\end{center}
\vskip -0.3in
\caption{Partial scan reconstructions of the Defrise phantom.} \label{fig:Defpartialscan}
\end{figure}

Half-scan reconstructions were done using data from only the eastern
hemisphere ($N_{\phi}=200$, $N_{\theta}=200$) or the southern
hemisphere ($N_{\phi}=400$, $N_{\theta}=100$). These hemispheres are
highlighted in Figure~\ref{fig:Defpartialscan}. The rest of the data
has been zero-filled.

It is known~\cite{Q1993,KLM,ref:LouisNQuinto,ref:Palamodov, ref:AnasPanTIP, ref:XWAK} that in case of incomplete data one can expect to recover stably only certain parts of the image the rest of it being blurred out. Namely some parts of the wavefront set of the image will be lost. For our phantom the singularities are jump discontinuities (edges) of imaged value $f$ across an interface $I$ (a surface of an ellipsoid in 3D). The wavefront $WF(f)$ of $f$ in this situation is the set of pairs $(x,n)$, where $x$ is a point on $I$, and $n$ is a vector normal to $I$ at $x$. As it was shown in~\cite{ref:LouisNQuinto, ref:XWAK} using microlocal analysis, a point $(x,n)\in WF(f)$ can be stably detected from the Radon data, if and only if $Rf$ includes data obtained from a sphere passing through $x$ and normal to $n$. In other words, one can see only those parts of an interface, that can be tangentially touched by spheres of integration centered at available transducer locations. The rest of the interface will be blurred.

Edges in the Defrise phantom were reconstructed in Figure~\ref{fig:Defpartialscan} as expected. When the data is collected from the eastern hemisphere there are enough spheres to touch tangentially all edges in the eastern hemisphere (see fig.~\ref{fig:Def2d}) but none to do it in the western hemisphere. That is why the locations of the edges in the eastern hemisphere were correctly reconstructed while those in the western part were blurred. When the data is collected from the southern hemisphere there are enough  spheres to touch tangentially all edges in the Defrise phantom, hence all of them were resolved.

From the geometric description above it is not hard to see that there may exist certain regions of reconstruction (locations of $x$, sometimes called audible zones) where any possible pair $(x,n)$ belonging to $WF(f)$ is recognizable from $Rf$. In our examples, when the data is collected from the eastern or southern hemisphere,  these regions are the eastern and southern half of the unit ball correspondingly.

Notice that the image values were not reconstructed correctly, since part of the data was missing. However, certain iterative techniques allow one to improve substantially the image values in the audible zone~\cite{ref:AnasPanTIP,ref:XWAK,ref:Paltauf}.


\subsection{Comparison with an approximate backprojection}

In early experimental work on thermoacoustic tomography, an approximate backprojection formula was used. It was written in analogy with the backprojection of regular Radon transform and looked similar to equation~(\ref{eq:rhofilt}), except the missing weight factor $\frac{1}{|\ex - \pea|}$.
The composition of this operator with the direct Radon transform is an elliptic pseudo-differential operator of order
zero (see, e.g.,~\cite{ref:LouisNQuinto, Guillemin, KLM}.)
Thus the locations and ``strengths" of image singularities are recovered correctly. However the values of the
image function will not be recovered correctly. The obtained reconstructions validate the predictions correctly recovering locations of edges. The values of image functions are accurate near the center where \mbox{$r\sim1$} but degrade slowly with distance from the origin, as expected.

\begin{figure}[h]
\begin{center}
\subfigure[Reconstruction of Defrise phantom at $256 \times 256$ resolution along the center \mbox{$x=0$} slice]{\epsfig{figure= 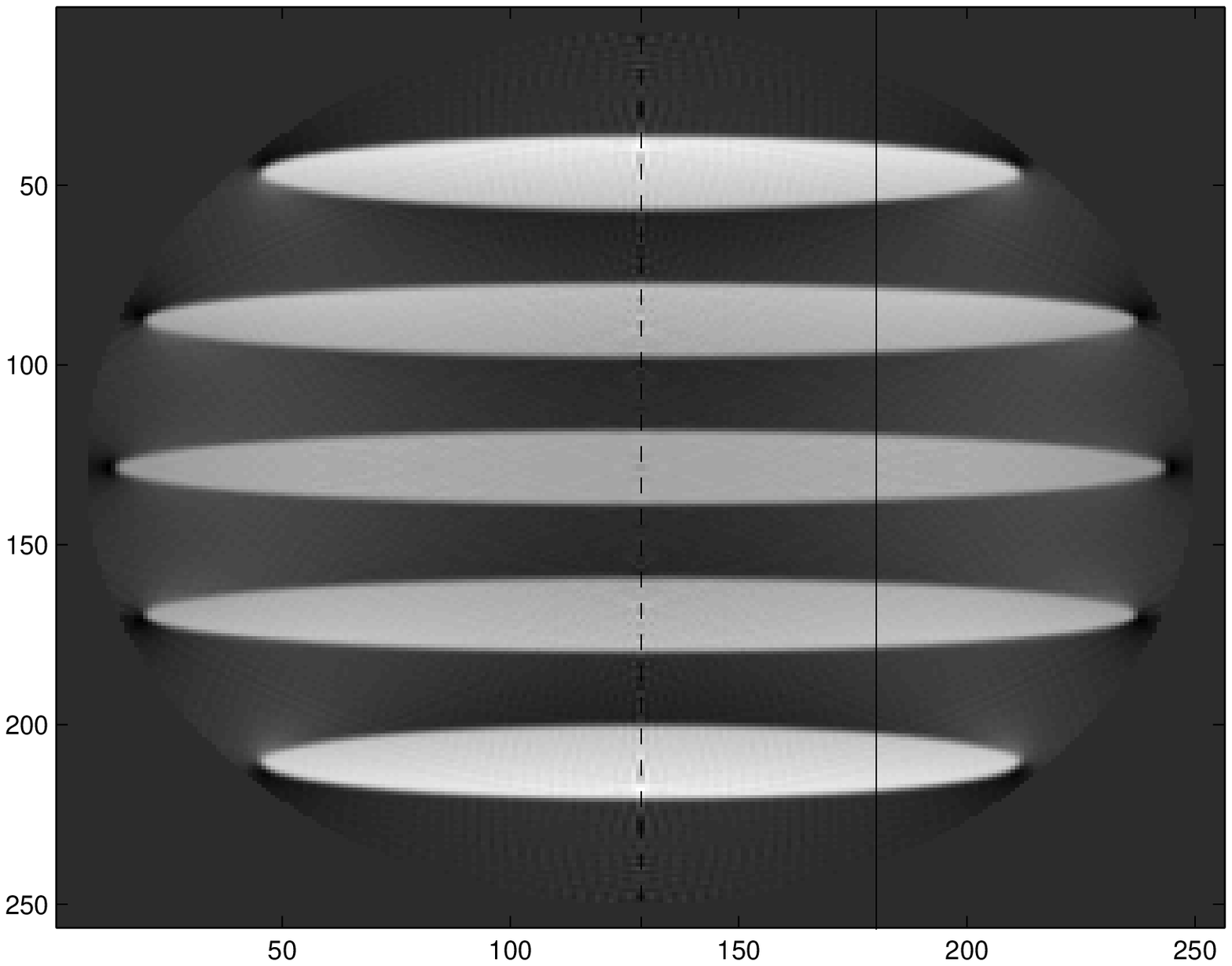,width=1.7in,height=1.7in }}
\subfigure[Sphere of radius 0.7 recon'd at low-res $64 \times 64$.  Below: center profiles \mbox{$y=0$} on three different slices.  Solid ~ \mbox{$x= \pm 0.1$}, dashed ~ \mbox{$x= \pm 0.3$}, dash-dot ~ \mbox{$x= \pm 0.5$}]{\epsfig{figure=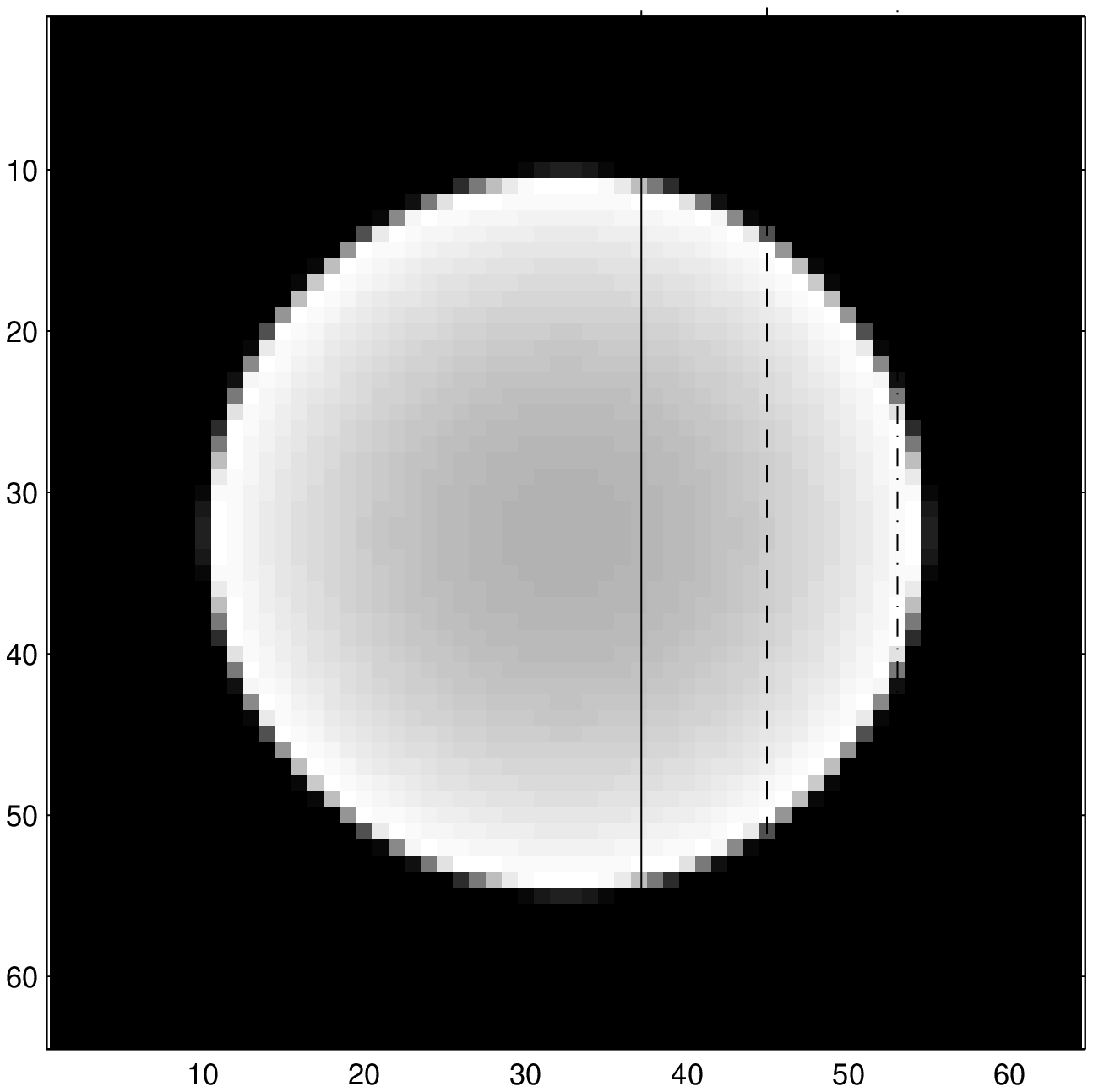,width=1.7in,height=1.7in }}
\\
{\epsfig{figure= 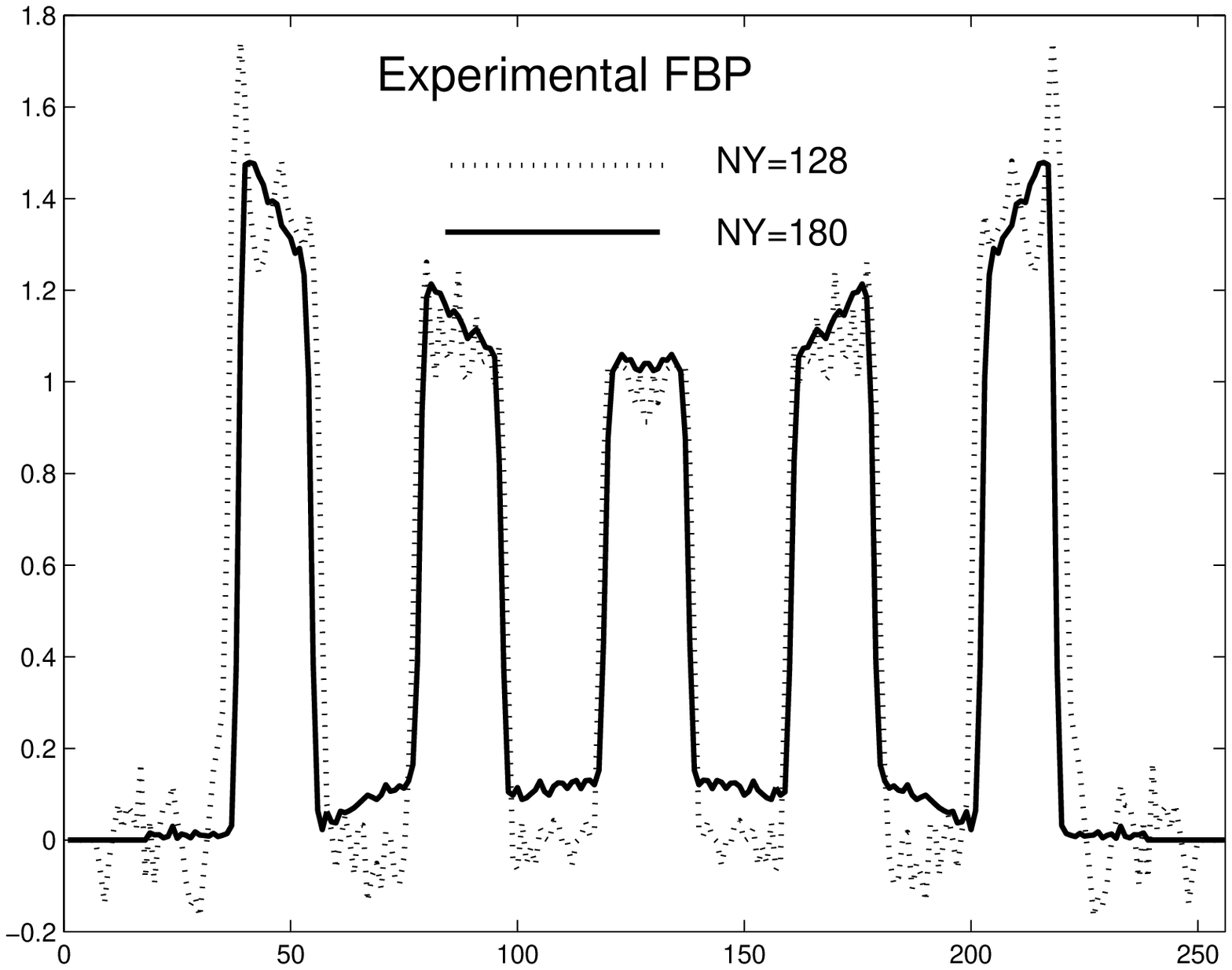,width=1.7in,height=1.7in }}
{\epsfig{figure= 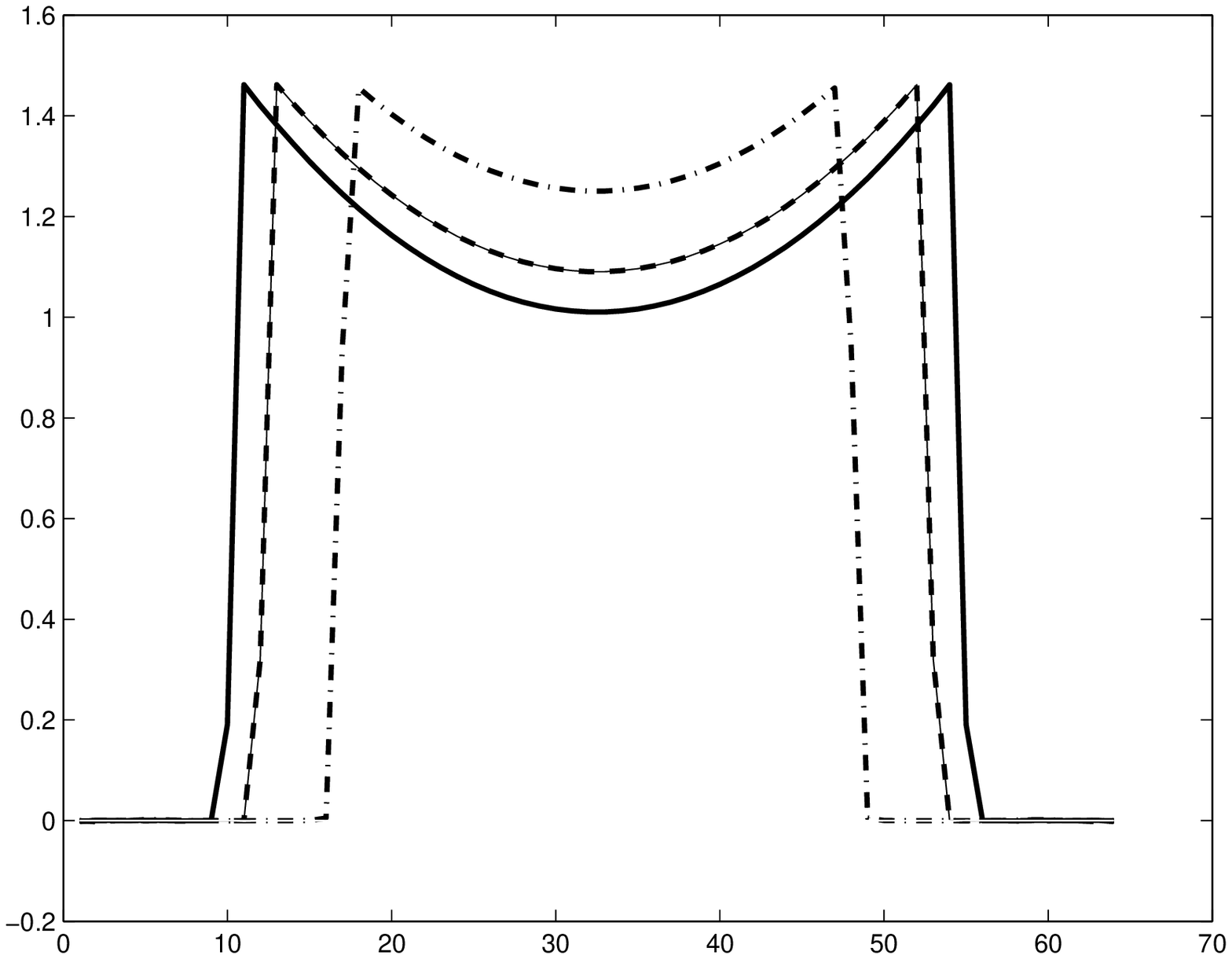,width=1.7in,height=1.7in }}
\end{center}
\vskip -0.3in
\caption{Approximate FBP shows low-frequency shading.}\label{exp-FBP}
\end{figure}


\subsection{Errors in reconstruction}\label{subsec:errors}

\begin{figure}[h]
\begin{center}
{\epsfig{figure= 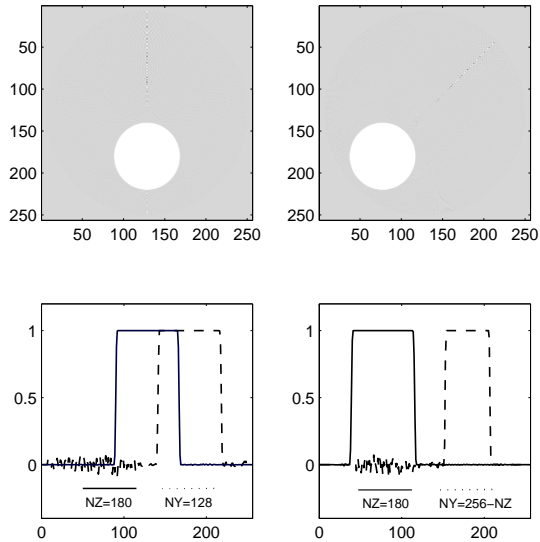,width=3.3in,height=3.3in }}
\end{center}
\vspace{-1.2cm} \caption{FBP errors along the axis of
symmetry.}\label{errors}
\end{figure}

As it was mentioned before, the reconstructions of Defrise phantom have some noise along the axis of phantom's symmetry $x=y=0$ (see figs.~\ref{Defrise}, \ref{fig:Defpartialscan}). To discuss the reasons of appearance of that noise we consider reconstructions of some simpler phantoms consisting of indicator functions of a perfect ball. This allows us to compute the Radon transform analytically, hence to exclude the errors in the data simulation. For every fixed $p_0$, $Rf(p_0,r)$ is a third order polynomial with respect to $r$ for $0<r_1\le r\le r_2<1$ and is zero for every other $r$. 
Filtered backprojection requires differentiating with respect to the radial variable $r$.  We used centered finite differences to estimate \mbox{$d^2/dr^2 Rf(p,r)$} which is exact on the third degree polynomials.  Therefore, we compute \mbox{$d^2/dr^2 Rf(p,r)$} exactly for all radii, $r$, at least $2\Delta r$ away from $r_1$ and $r_2$. Hence the only errors in numerical differentiation that spread into the backprojection come with the data from spheres close to the ones touching tangentially the phantom ball. None of these spheres passes inside the phantom ball, hence backprojection at those points is free of errors from numerical differentiation (see fig.~\ref{errors}).

\begin{figure}[h]
\begin{center}
{\epsfig{figure= 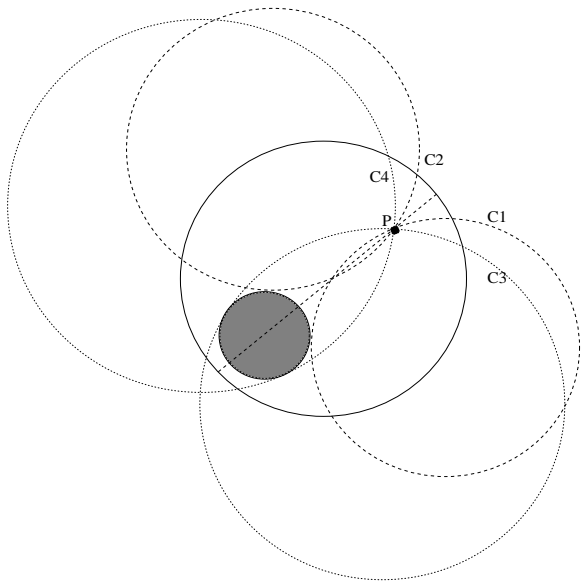,width=3.0in,height=3.0in }}
\end{center}
\vspace{-0.5cm}
\caption{The spheres C3 and C4 contain the phantom ball, while C1 and C2 do not}\label{circles}
\end{figure}

Now let us consider a point $p_1$ on the axis of symmetry of the ball phantoms (the line connecting the center of the phantom ball and the origin). There are two sets of spheres that pass through that point and touch the phantom ball tangentially. The spheres in the first set contain the phantom ball, while the spheres in the second set do not.
A 2D slice of this scenario is presented in Figure~\ref{circles}. 

Notice that all spheres in the same set have the same radius. So the
errors from the numerical differentiation that they will bring into
the backprojection algorithm are absolutely the same. The axis of
symmetry is the only location in the reconstruction region where
these errors are perfectly correlated. This resonance increases the
magnitude of errors resulting in the noise along the symmetry axis
on reconstructed images (see fig.~\ref{errors}).

\begin{figure}[h]
\begin{center}
{\epsfig{figure= 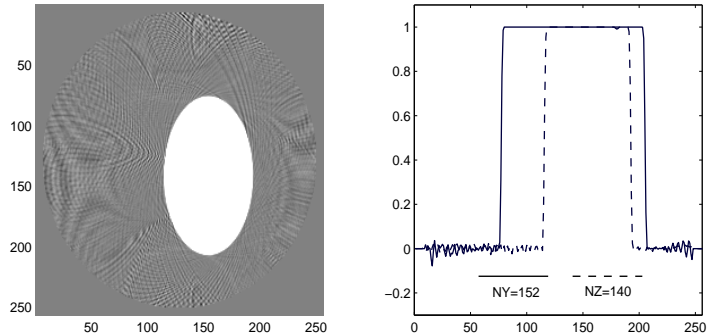,width=4.5in,height=2in }}
\end{center}
\vspace{-0.5cm}
\caption{An ellipsoidal phantom with center at (0,0.2,-0.1) and semiaxes lengths equal to (0.4,0.3,0.5).}\label{ellipsoid}
\end{figure}

In case of ellipsoids in the Defrise phantom everything said above holds. In fact magnitude of errors is five times bigger since there are five ellipsoids with the same axis of symmetry there. At the same time, the reconstruction of an ellipsoidal phantom without any rotational symmetry has no axis of emphasized errors (see fig.~\ref{ellipsoid}).

\section{Conclusion}

We have implemented a straightforward numerical validation of both
FBP and $\rho$-filtered inversion formulae for TCT data on the
high-frequency Defrise phantom. FBP and $\rho$-filtered have
virtually identical performance with noise-free simulated data from
this high-contrast object.  Artifacts due to numerical errors are
more severe in FBP than $\rho$-filtered images and might be reduced
by mollification techniques~\cite{ref:HaltSchuScherzer}.  Comparing
FBP and $\rho$-filtered performance in the presence of noise and for
low-contrast detectability will appear in future publications.


\section{Acknowledgements}

The authors would like to thank Mark Anastasio, David Finch, Peter
Kuchment, Leonid Kunyansky and Rakesh for information about their
work and discussions on the subject. The first author would also
like to thank his advisor Peter Kuchment for constant support.

The work has been done in Summer 2004 at GE Healthcare Technologies,
Milwaukee, WI as part of the internship program of the first author.

The first author was supported in part by the NSF Grants DMS 9971674
and 0002195 and thanks the NSF for that support. Any opinions,
findings, conclusions or recommendations expressed in this paper
are those of the authors and do not necessarily reflect the views of
the NSF.

\end{document}